\documentclass[12pt]{article}

\usepackage{epsf,graphicx}
\usepackage{cite, amsmath, amssymb}
\usepackage{tikz,color}
\newtheorem{theorem}{\bf Theorem}[section]
\newtheorem{corollary}[theorem]{\bf Corollary}
\newtheorem{lemma}[theorem]{\bf Lemma}

\newcommand{\proof}{\noindent{\bf Proof.\ }}
\newcommand{\qed}{\hfill $\square$ \bigskip}
\newcommand{\ecc}{{\rm ecc}}

\newcommand{\diam}{{\rm diam}}
\newcommand{\Tr}{{\rm Tr}}
\newcommand{\Mo}{{\rm Mo}}
\newcommand{\Z}{{\rm Z}}

\begin{document}

\title{Extremal results on stepwise transmission irregular graphs}

\author{Yaser Alizadeh$^{a}$  \and Sandi Klav\v zar$^{b,c,d}$
}

\date{\today}

\maketitle
% \vspace{-0.8 cm}
\begin{center}

$^a$ Department of Mathematics, Hakim Sabzevari University, Sabzevar, Iran \\
e-mail: {\tt y.alizadeh@hsu.ac.ir} \\
\medskip
 
$^b$ Faculty of Mathematics and Physics, University of Ljubljana, Slovenia\\
e-mail: {\tt sandi.klavzar@fmf.uni-lj.si} \\
\medskip

$^c$ Faculty of Natural Sciences and Mathematics, University of Maribor, Slovenia \\
\medskip

$^d$ Institute of Mathematics, Physics and Mechanics, Ljubljana, Slovenia\\
\end{center}

\begin{abstract} 
The transmission ${\rm Tr}_G(v)$ of a vertex $v$ of a connected graph $G$ is the sum of distances between $v$ and all other vertices in $G$. $G$ is a stepwise transmission irregular (STI) graph if $|{\rm Tr}_G(u) - {\rm Tr}_G(v)| =1$ holds for each edge $uv \in E(G)$. In this paper, extremal results on STI graphs with respect to the size and different metric properties are proved. Two extremal families appear in all the cases, balanced complete bipartite graphs of odd order and the so called odd hatted cycles. 
\end{abstract}


\medskip\noindent
{\bf AMS Subj.\ Class. 2020:} 05C12 

%%%%%%%%%%%%%%%%%%%%%%%%%%%%%%%%%%%%%%%
\section{Introduction}
\label{sec:intro}
%%%%%%%%%%%%%%%%%%%%%%%%%%%%%%%%%%%%%%%

All graphs throughout the paper are simple and connected.  If $G = (V(G), E(G))$ is a graph and $u,v\in V(G)$, then the distance $d_G(u,v)$  is the number of edges on a shortest $u,v$-path. The {\em transmission} $\Tr_G(v)$ of a vertex $v\in V(G)$ is the sum of distances between $v$ and all other vertices. This concept arose in different contexts, hence it is not surprising that it is also known as the total distance, the farness, and the vertex Wiener value, cf.~\cite{abiad-2017, hua-2018, krnc-2015, rajasingh-2016, sharafdini-2020}. The transmission also led to the Wiener complexity (the number of different transmissions)~\cite{alizadeh-2014}, is closely related to other topological indices~\cite{sharafdini-2020}, and characterizes the distance-balanced property and the opportunity index~\cite{balakrishnan-2014, cavaleri-2019}. 

Quite recently, Dobrynin and Sharafdini~\cite{dob1} defined {\em stepwise transmission irregular graphs} ({\em STI graphs} for short) as the graphs $G$ in which $|\Tr_G(u)-\Tr_G(v)|=1$ holds for each edge $uv$ of $G$. The class was proposed, among other things, as a new pebble in the investigation of transmission dependent classes of graphs. A striking example of such classes is the one of transmission irregular graphs which by definition contains the graphs in which all its vertices have pairwise different transmissions. After being introduced in~\cite{AK2018}, transmission irregular graphs received a lot of attention~\cite{al-yakoob-2020, al-yakoob-2021, bezhaev-2021, dobrynin-2019-b, dobrynin-2019-c, xu-2021a}. 

Just as transmission irregular graphs, STI graphs also turned out to be a very interesting class of graphs. In the seminal paper~\cite{dob1}, basic properties of STI graphs were first established. In particular, STI graphs are bipartite, $2$-connected (except $P_3$), and of odd order, so that no regular graph can be STI. By computer search it was established that there are $1$, $1$, $3$, $7$, $18$, $87$, and $1171$ STI graphs of respective orders  $3$, $5$, $7$, $9$, $11$, $13$, and $15$. Six interesting infinite families of STI graphs were also presented and topological indices of STI graphs briefly discussed. In the subsequent paper~\cite{al-yakoob-2022},  Al-Yakoob and Stevanovi\'c confirmed a conjecture from~\cite{dob1} by proving that all the graphs from an additional interesting infinite family are STI. They also provided a further computational support for a challenging conjecture from~\cite{dob1} that the girth of every STI graph is $4$. As feasible candidates for a counterexample they checked all bipartite, $2$-connected graphs on up to $23$ vertices with girth at least $6$.  

In this paper we are interested in extremal properties of STI graphs. In the second part of this section definitions and notation needed is given. In Section~\ref{sec:size} we determine the minimum and the maximum size of an STI graph and characterize the extremal graphs. In Section~\ref{sec:metric} we determine extremal results for different metric properties of STI graphs: the diameter, the eccentricity, the Wiener index, and the transmission. 

%%%%%%%%%%%%%%%%%%%%%%%%%%%%%%%%%%%%%%%
\subsection{Preliminaries}
%%%%%%%%%%%%%%%%%%%%%%%%%%%%%%%%%%%%%%%

The {\em eccentricity} $\varepsilon_G(v)$ {\em of a vertex} $v\in V(G)$  is the maximum distance between $v$ and other vertices of $G$.  The {\em eccentricity}  $\ecc(G)$ {\em of a graph} $G$ is the sum of the vertex eccentricities over all vertices of $G$. A vertex $w$ is an {\em eccentric vertex of} $v$ if $d_G(v,w) = \varepsilon_G(v)$. The {\em diameter} $\diam(G)$ of $G$ is the largest eccentricity among its vertices. 

The {\em Wiener index} $W(G)$ of a graph $G$ is the sum of distances between all pairs vertices of $G$, that is, $W(G) = \sum_{\{u,v\}\subseteq V(G)}d_G(u,v)= \frac{1}{2}\sum_{v\in V(G)}\Tr_G(v)$, see the  survey~\cite{knor-2016}. For an edge $vw\in E(G)$, the number of vertices lying closer to $v$ than to $w$ is denoted by $n_v(G)$. The value $n_w(G)$ is defined analogously. The {\em Mostar index} $\Mo(G)$ of $G$ is defined as $\Mo(G)= \sum_{uv \in E(G)}|n_u(G) -n_v(G)|$. It was introduced in~\cite{Tom2} and studied a lot afterwards, see the recent survey~\cite{ali-2021}. The {\em first Zagreb index} $\Z_1(G)$ of $G$ is defined as $\Z_1(G) = \sum_{u\in V(G)}\deg(u)^2$, the {\em second Zagreb index} $\Z_2(G)$ of $G$ is $\Z_2(G) = \sum_{uv\in E(G)}\deg(u)\deg(v)$. These two invariants were respectively introduced in~\cite{Gut1, Gut2}, see also the survey~\cite{gutman-2013}. We note that the first Zagreb index can be equivalently expressed as $\Z_1(G) = \sum_{uv\in E(G)}\deg(u)+\deg(v)$, see~\cite{Tom}.

For further reference we recall the following properties of STI graphs. 

\begin{theorem} \label{basic}{\rm \cite{dob1}}
If $G$ is an STI graph different from $P_3$, then $G$ is a bipartite, $2$-connected graph of odd order.   
\end{theorem}

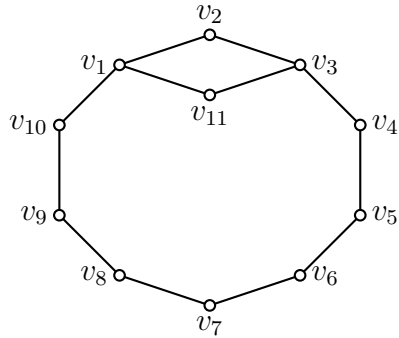
\begin{figure}[ht!]
\begin{center}
\begin{tikzpicture}[scale=0.8,style=thick]
\tikzstyle{every node}=[draw=none,fill=none]
\def\vr{2.5pt} % \vr = vertex radius;  Set \vr = 2/scale for uniform sizing of vertices

\begin{scope}[yshift = 0cm, xshift = 0cm]
%% vertices defined %%
\path (1,4) coordinate (v1);
\path (2.5,4.5) coordinate (v2);
\path (4,4) coordinate (v3);
\path (5,3) coordinate (v4);
\path (5,1.5) coordinate (v5);
\path (4,0.5) coordinate (v6);
\path (2.5,0) coordinate (v7);
\path (1,0.5) coordinate (v8);
\path (0,1.5) coordinate (v9);
\path (0,3) coordinate (v10);
\path (2.5,3.5) coordinate (v11);

%% edges %%
\draw (v1) -- (v2) -- (v3) -- (v4) -- (v5) -- (v6) -- (v7) -- (v8) -- (v9) -- (v10) -- (v1) -- (v11) -- (v3);
%% vertices %%%
\draw (v1)  [fill=white] circle (\vr);
\draw (v2)  [fill=white] circle (\vr);
\draw (v3)  [fill=white] circle (\vr);
\draw (v4)  [fill=white] circle (\vr);
\draw (v5)  [fill=white] circle (\vr);
\draw (v6)  [fill=white] circle (\vr);
\draw (v7)  [fill=white] circle (\vr);
\draw (v8)  [fill=white] circle (\vr);
\draw (v9)  [fill=white] circle (\vr);
\draw (v10)  [fill=white] circle (\vr);
\draw (v11)  [fill=white] circle (\vr);
%% text %%
\draw[left] (v1) node {$v_1$};
\draw[above] (v2) node {$v_2$};
\draw[right] (v3) node {$v_3$};
\draw[right] (v4) node {$v_4$};
\draw[right] (v5) node {$v_5$};
\draw[right] (v6) node {$v_6$};
\draw[below] (v7) node {$v_7$};
\draw[left] (v8) node {$v_8$};
\draw[left] (v9) node {$v_9$};
\draw[left] (v10) node {$v_{10}$};
\draw[below] (v11) node {$v_{11}$};
\end{scope}

\end{tikzpicture}
\end{center}
\caption{The graph $G_{11}$}
\label{fig1}
\end{figure}

To conclude the preliminaries we introduce a family of graphs that will be ubiquitous throughout the rest of the paper. Let $n\ge 5$. The graph $G_n$ has the vertex set $\{v_1, \ldots, v_n\}$, where the vertices $v_1, \ldots, v_{n-1}$ induce a cycle of length $n-1$, while $v_n$ is adjacent to $v_{1}$ and $v_3$. Note that $G_5 \cong K_{2,3}$, while for $G_{11}$ see Fig.~\ref{fig1}. 

We note in passing that in~\cite{bujtas-2021}, the graphs $G_n$ were named {\em hatted cycles} and investigated with respect to the Rall's $1/2$-conjecture on the domination game. From our current point of view, it was shown in~\cite{dob1} that $G_n$ is an STI graph for every odd $n \ge 5$. Moreover, in~\cite{bessy-2020} it was proved that hatted cycles are the graphs that attain the third maximum Wiener index over $2$-connected graphs.

%%%%%%%%%%%%%%%%%%%%%%%%%%%%%%%%%%%%%%%
\section{Extremal STI graph with respect to size}
\label{sec:size}
%%%%%%%%%%%%%%%%%%%%%%%%%%%%%%%%%%%%%%%

In this section we bound the size of an STI graph and determine the extremal graphs. As a consequence we do the same for the Mostar index and the two Zagreb indices.

\begin{theorem}\label{edge}
If $G$ is an STI graph of order $n$ and size $m$, then 
$$n + 1  \le  m  \le \frac{n^2-1}{4}\,.$$
Moreover, the left equality holds if and only if $G \cong G_n$, and the right equality holds if and only if $G \cong  K_{\frac{n-1}{2},\frac{n+1}{2}}$.
\end{theorem}

\proof 
The right inequality was proved in~\cite[Proposition 4]{dob1}. For the equality, let $G$ be an arbitrary STI graph of order $n$ and size $m$. Then, by Theorem~\ref{basic}, $G$ is bipartite and $n$ is odd. Since $K_{\frac{n-1}{2},\frac{n+1}{2}}$ is the unique bipartite graph of order $n$ with the maximum size and $K_{\frac{n-1}{2},\frac{n+1}{2}}$ is an STI graph (cf.~\cite{dob1}), it is thus the unique graph attaining the right equality.  

Recall that a graph $G$ is $2$-connected if and only if $G$ admits an ear decomposition  (see~\cite{west-2001} for the definition of the ear decomposition). Hence the only $2$-connected graph $G$ with $m = n$ is $C_n$. As cycles are not STI graphs, we conclude that if $G$ is an STI graph, then $m  \ge n + 1$. Let now $G$ be an arbitrary STI graph with $m = n + 1$. Then it is $2$-connected and the ear decomposition implies that $G$ consists of two vertices of degree $3$, say $x$ and $y$, which are connected by three internally disjoint paths. Let these paths be of lengths $r,s,t$, where $1\le r \le s \le t$. 

Since $G$ is bipartite, $r$, $s$, and $t$ are either all odd or all even. Suppose first that $r=1$, that is, $x$ and $y$ are connected by an edge. Then clearly $\Tr_G(x) = \Tr_G(y)$ and so $G$ is not an STI graph.  
Assume in the rest that $r \ge 2$. If $r=s=2$, then $G\cong G_n$ which was shown to be an STI graph in~\cite[Proposition~7]{dob1}. Hence let $s\ge 3$. If $P_1$, $P_2$, and $P_3$ are the $x,y$-paths of respective lengths $r$, $s$, and $t$, and $w$ is the neighbor of $x$ on $P_3$, then we can estimate as follows: 
\begin{eqnarray*}
\Tr_G(w) - \Tr_G(x) &=& \sum_{z\in{P_1}}(d_G(w,z) - d_G(x,z))+\sum_{z\in{P_2\cup P_3}\setminus \{x,y\}}(d_G(w,z) - d_G(x,z)) \\ &\ge & \sum_{z\in{P_1}}(d_G(w,z) - d_G(x,z)) \ge 2\,.
\end{eqnarray*}
We conclude that $G_n$ is the only STI graph with $m = n + 1$.
\qed

\begin{corollary}
If $G$ is an STI graph of order $n\ge 5$, then the following holds. 
\begin{enumerate}
\item[(i)] $\displaystyle{  n+1 \le \Mo(G) \le \frac{n^2-1}{4}\,.}$
\item[(ii)]  $\displaystyle{ 4n + 1 \le \Z_1(G) \le \frac{n(n^2-1)}{4}\,.}$
\item[(iii)] $\displaystyle{ 4n+16 \le \Z_2(G) \le \left(\frac{n^2-1}{4}\right)^2\,. }$
\end{enumerate}
Moreover, in each of the cases the left equality holds if and only if $G=G_n$ and the right equality holds if and only if $G=K_{\frac{n-1}{2},\frac{n+1}{2}}$.
\end{corollary}

\proof
(i) This follows from the fact that if $G$ is an arbitrary graph, then $\Mo(G)$ can be equivalently written as $\sum_{uv \in E(G)}|\Tr_G(u)-\Tr_G(v)|$, see~\cite[Corollary 2.2]{aliz}.

(ii) By Theorem~\ref{edge}, and since $G$ is $2$-connected, $G$ contains at least two vertices of degree $3$, while the other vertices are of degree at least $2$. Thus 
\[\Z_1(G) = \sum_{v \in V(G)}\deg(v)^2 \ge 4(n-2)+9 = 4n+1\,.\]
Equality holds here if and only if $G$ has exactly two vertices of degree $3$ and the other vertices are of degree $2$ which means (having the ear decomposition in mind) that $G\cong G_n$. For the right inequality observe first that if $uv$ is an edge of (a bipartite graph) $G$, then $\deg(u) + \deg(v) \le n$. Therefore,
\[\Z_1(G)= \sum_{uv\in E(G)}\deg(u) + \deg(v) \le n |E(G)|\le \frac{n(n^2-1)}{4}.\] 
Using Theorem~\ref{edge} again we infer that the right equality holds if and only if $G \cong K_{\frac{n-1}{2},\frac{n+1}{2}}$.

(iii) This is proved using a similar argument as (ii).
\qed

%%%%%%%%%%%%%%%%%%%%%%%%%%%%%%%%%%%%%%%
\section{Metric properties}
\label{sec:metric}
%%%%%%%%%%%%%%%%%%%%%%%%%%%%%%%%%%%%%%%

In this section we consider extremal behavior of STI graphs with respect to metric invariants. In the first subsection we begin with their diameter and determine the extremal graphs. Based on this result, the extremal STI graphs are determined with respect to the eccentricity. In the second subsection we prove parallel results for the Wiener index and the transmission.  

%%%%%%%%%%%%%%%%%%%%%%%%%%%%%%%%%%%%%%%
\subsection{Eccentricity}
%%%%%%%%%%%%%%%%%%%%%%%%%%%%%%%%%%%%%%%

\begin{lemma}\label{diameter}
If $G$ is an STI graph of order $n$, then
\[ 2\le \diam(G) \le \frac{n-1}{2}\,.\]
Moreover, the left equality holds if and only if $G \cong K_{\frac{n-1}{2},\frac{n+1}{2}}$, and the right equality holds for $G_n$. 
\end{lemma}

\proof
The left inequality follows from the fact that the only bipartite graph of diameter $1$ is $K_2$ which is not an STI graph. If the left equality holds, that is, if $\diam(G) = 2$, then $G$ is a complete bipartite graph. And since $G$ is an STI graph, we conclude that $G \cong K_{\frac{n-1}{2},\frac{n+1}{2}}$. 

By Theorem~\ref{basic}, $G$ is a 2-connected graph of odd order which in turn implies the right inequality. It is clear that $\diam(G_n) = (n-1)/2$.  \qed

\begin{theorem}\label{eccentricity}
If $G$ is an STI graph of order $n\ge 4$, then
$$ 2n \le \ecc(G) \le \frac{n(n-1)}{2}\,.$$
Moreover, the left equality holds if and only if $G \cong K_{\frac{n-1}{2},\frac{n+1}{2}}$ and the right equality holds if and only if $G \cong G_n$. 
\end{theorem}

\proof
By Lemma~\ref{diameter}, $2\le \diam(G) \le (n-1)/2$. 

If $\diam(G) = 2$, then Theorem~\ref{basic} implies that $G$ is a complete bipartite graph. Among them only $K_{\frac{n-1}{2},\frac{n+1}{2}}$ is an STI graph.

As already said, $\diam(G) \le (n-1)/2$ which in turn implies that $\ecc(G) \le n \frac{(n-1)}{2}$. Suppose now that the right equality holds for an STI graph $G$. Then each vertex $v$ has $\varepsilon_G(v) = (n-1)/2$. Let $v$ and $w$ be two vertices with $d_G(v,w) = (n-1)/2$. As $G$ is 2-connected, there exists internally disjoint $v,w$-paths $P_1$ and $P_2$. Then each of these paths is of length at least $(n-1)/2$, which means that $|V(P_1)\cup V(P_2)| \ge n-1$. We now distinguish two cases. 

Suppose first that $|V(P_1)|= (n-1)/2$ and $|V(P_2)|= (n+1)/2$. In this case, $V(G) = V(P_1) \cup V(P_2)$. Since $G$ is an STI graph and cycles are not, there exists an edge $pq \notin E(P_1\cup P_2)$. As $G$ is bipartite, $d_{P_1\cup P_2}(p,q)\ge 3$. But this implies that $\varepsilon_G(p) < (n-1)/2$ which is not possible. In this case we thus have no equality situation. 

Suppose second that $|V(P_1)|=|V(P_2)| = (n-1)/2$. Then there exists exactly one vertex $x \notin V(P_1\cup P_2)$. As $G$ is $2$-connected, $x$ has at least two neighbors in $V(P_1)\cup V(P_2)$, let $y$ and $z$ be its arbitrary neighbors. If $d_{P_1\cup P_2}(y,z) \ge 3$, then $\varepsilon_G(y) < (n-1)/2$. Hence we must have $d_{P_1\cup P_2}(y,z)=2$. If $\deg(x)=3$, then $\varepsilon_G(y)< (n-1)/2$ or $\varepsilon_G(z)< (n-1)/2$. Also, if there exists an edge not in $E(P_1) \cup E(P_2)\cup\{xy, xz\}$, say $st$, then $\varepsilon_G(s) < (n-1)/2$. We conclude that $G \cong G_n$. 
\qed

%%%%%%%%%%%%%%%%%%%%%%%%%%%%%%%%%%%%%%%
\subsection{Wiener index and transmission}
%%%%%%%%%%%%%%%%%%%%%%%%%%%%%%%%%%%%%%%

\begin{theorem}\label{wiener}
If $G$ is an STI graph of order $n$, then 
$$\frac{3n^2-4n+1}{4} \le W(G) \le \frac{n^3 -n^2-n+17}{8}\,.$$
Moreover, the left equality holds if and only if $G \cong K_{\frac{n-1}{2},\frac{n+1}{2}}$, and the right equality holds if and only if $G \cong G_n$.
\end{theorem}

\proof
By Theorem~\ref{basic}, $G$ is bipartite. Let $V(G) = V_1 \cup V_2$ be its bipartition. Since adding an edge between two nonadjacent vertices decrease the Wiener index of the graph, we have $W(G) \ge W(K_{|V_1|,|V_2|})$. Moreover, $K_{\frac{n-1}{2},\frac{n+1}{2}}$ is the only complete bipartite of order $n$ which is also an STI graph, hence the left inequality holds, as well as its equality part.  

For the right inequality, let $H_{n,p,q}$ be a graph of order $n$ comprised of three internally disjoint paths with the same end-vertices and of respective lengths $p$, $q$, and $n - p - q + 1$. In~\cite[Theorem 1]{bessy-2020} it was proved that if $n\ge 11$, and $G$ is a $2$-connected graph of order $n$ different from $C_n$, $H_{n,1,2}$, and $H_{n,2,2}$, then $W(G) < W(H_{n,2,2}) < W(H_{n,1,2}) < W(C_n)$. Note now that $G_n \cong H_{n,2,2}$ and that $H_{n,1,2}$ is not an STI graph. Hence, as for odd $n$ we have $W(G_n) = (n^3 -n^2-n+17)/8$, cf.~\cite[Table 2]{bessy-2020}, the right inequality follows together with the equality part. For the latter, one needs to consider the small cases separately, that is, \cite[Theorem 7]{bessy-2020} for the case $n=9$ and seven sporadic STI graphs with $n=7$ from~\cite[Figure 1]{dob1}.   
\qed  
  
\begin{theorem}\label{max-transmission}
If $G$ is an STI graph of order $n$, then 
\[\max \{\Tr_G(v):\ v \in V(G)\} \le \frac{n^2-1}{4}\,.\]
Moreover, the equality holds if and only if $G \cong G_n$.
\end{theorem}

\proof
Let $v$ be a vertex of $G$ and consider the BFS-tree rooted at $v$. Since $G$ is 2-connected, there must be at least two vertices in each of the levels of the BFS-tree, that is, there exists at least two vertices at each possible distance from $v$. Since $n$ is odd, this implies that 
$$\Tr_G(v) \le 2(1+\cdots + \frac{n-1}{2}) = \frac{n^2-1}{4}\,,$$
which proves the upper bound.
  
The equality holds if and only if for each positive integer $k$, $1\le k \le (n-1)/2$, there are exactly two vertices at distance $k$ from $v$ and $\varepsilon_G(v) = (n-1)/2$. Let $G$ be an STI graph for which the equality holds, and let $v\in V(G)$ be a vertex with $\Tr_G(v) = (n^2-1)/4$.  Let $w$ be an eccentric vertex of $v$, so that $d_G(v,w) = (n-1)/2$. As $G$ is $2$-connected, there exist two internally vertex disjoint $v,w$-paths $P_1$ and $P_2$. Then $V(G) = V(P_1) \cup V(P_2) \cup \{w'\}$, where $w'$ is the second eccentric vertex of $v$. Using the 2-connectivity of $G$ again, we infer that $w'$ is adjacent to the two vertices at distance $(n-3)/2$ from $v$. Now, if these two edges, together with the edges of $P_1$ and $P_2$ are all the edges of $G$, then $G\cong G_n$. Suppose next that $G$ contains some other edge $f$. Since $P_1$ and $P_2$ are shortest paths, such an edge connects a vertex $x$ of $P_1$ with a vertex $y$ of $P_2$. Suppose that $x$ is selected such that it is closest to $v$ among all veretices that are endpoints of such additional edges $f$. As $G$ is bipartite this implies that $d_G(v,y) = d_G(v,x) + 1$. Let $x'$ be the predecessor of $x$ on $P_1$ (so that $d_G(v,x) = d_G(v,x') + 1$). It is possible that $x' = v$. 

We claim that $\Tr(x') - \Tr(x)\ge 2$. By the way the edge $xy$ is selected, the vertices that are at distance at most $d_G(v,x)$ from $v$ together with the vertex $y$ induce a cycle $C$.  Moreover, again by the way the edge $xy$ is selected, $C$ is an isometric cycle, that is, if $z,z'\in V(C)$, then $d_C(z,z') = d_G(z,z')$. It follows that the vertices of $C$ contribute the same value to $\Tr_G(x)$ and to $\Tr_G(x')$. Since for every vertex $z''$ from $V(G)\setminus V(C)$ we have $d_G(x,z'') < d_G(x',z'')$, the claim follows. 

From the above claim we conclude that there is no such edge $f$ in $G$ and hence $G\cong G_n$. To complete the argument observe that $G_n$ contains exactly one vertex, say $v$, with $\varepsilon_G(v) = (n-1)/2$, and for $v$ we have $\Tr_{G_n}(v) = (n^2-1)/4$. 
\qed

\begin{theorem}\label{min-transmission}
If $G$ is an STI graph of order $n$, then 
\[\min \{\Tr_G(v):\ v \in V(G)\} \ge \frac{3n-5}{2}\,.\]
Moreover, equality holds if and only if $G \cong K_{\frac{n-1}{2},\frac{n+1}{2}}$.
\end{theorem}

\proof 
Let $(X,Y)$ be the bipartition of $G$ and assume without loss of generality that $|X|\le |Y|$. Since the order $n$ of an STI graph $G$ is odd, we actually have $|X| < |Y|$. Let $v \in X$ has the minimum transmission among vertices of $G$. Then $\Tr_G(v) \ge |Y| + 2(|X|-1)$, and equality holds if and only if $v$ is adjacent to all vertices in $Y$ and is at distance $2$ from each vertex of $X$. Let $w \in Y$ be adjacent to $v$. Since $G$ is an STI graph and $|X| < |Y|$, we have $\Tr_G(w)=\Tr_G(v)+1 = |Y| + 2|X|-1$. If there would be a vertex $z \in X$ non adjacent to $w$, or a vertex $z \in Y$ with $d_G(w,z)\ge 4$, then 
\[\Tr_G(w) \ge (|X|-1) + 2(|Y|-1) + 3 = |X| + 2|Y| > \Tr_G(v) + 1\,,\]
which is not possible. Hence $w$ must be adjacent to all vertices of $X$ and at distance $2$ from each vertex in $Y$. As $w$ was an arbitrary neighbor of $v$ we conclude that $v$ has the minimum transmission $\Tr_G(v) = |Y| + 2(|X|-1)$ if and only if $G$ is a complete bipartite graph. We conclude that  $G \cong K_{\frac{n-1}{2},\frac{n+1}{2}}$ and $\Tr_G(v) = (3n-5)/2$.
\qed

\section*{Acknowledgements}

Sandi Klav\v{z}ar acknowledges the financial support from the Slovenian Research Agency (research core funding P1-0297 and projects N1-0095, J1-1693, J1-2452).

\end{document}